\documentclass[12pt,oneside,reqno]{amsart}

\usepackage[top=2.5cm,bottom=2.5cm,left=3cm,right=3cm]{geometry}

\usepackage{mathtools}
\usepackage{needspace}
\usepackage{amssymb} % Symbols used in the manuscript.
\IfFileExists{newtx.sty}{
  \usepackage[libertinus]{newtx}
}{
  \IfFileExists{libertinus.sty}{\usepackage{libertinus}}{
    \usepackage[T1]{fontenc}
    \usepackage{lmodern}
  }
}
\usepackage{microtype}
\usepackage{mathrsfs}
\DeclareFontFamily{U}{rsfs}{\skewchar\font127 }
\DeclareFontShape{U}{rsfs}{m}{n}{%
  <-6.5> rsfs5
  <6.5-8> rsfs7
  <8-> rsfs10
}{}
\usepackage[mathcal]{eucal}

\usepackage{tikz-cd}
\usepackage{enumitem}
\setlist{itemsep=1pt,topsep=4pt}
\usepackage{graphicx}
\usepackage{comment}
\DeclareMathAlphabet{\mathcalCM}{OMS}{cmsy}{m}{n}

\theoremstyle{plain}
\newtheorem{theorem}{Theorem}[section]
\newtheorem{proposition}[theorem]{Proposition}
\newtheorem{lemma}[theorem]{Lemma}
\newtheorem{corollary}[theorem]{Corollary}

\theoremstyle{definition}
\newtheorem{definition}[theorem]{Definition}

\theoremstyle{remark}
\newtheorem{remark}[theorem]{Remark}

\numberwithin{equation}{section}

\usepackage{xcolor}
\definecolor{citeblue}{rgb}{0,0.35,0.75}
\usepackage[
  colorlinks=true,
  linkcolor=citeblue,
  citecolor=citeblue,
  urlcolor=citeblue
]{hyperref}
\usepackage[nameinlink,noabbrev]{cleveref}
\Crefname{theorem}{Theorem}{Theorems}
\Crefname{conjecture}{Conjecture}{Conjectures}
\Crefname{question}{Question}{Questions}

\newcommand{\cF}{\mathcal F}
\newcommand{\PP}{\mathbb P}

\title{Chow Stability of Polarized K3 Surfaces with ADE singularities}
\author[Z. Li]{Zhiyuan Li}\address{Shanghai Center for Mathematical Sciences, Fudan University, 2005 Songhu Road 200438, Shanghai, China}\email{zhiyuan\_li@fudan.edu.cn}
\date{\today}

\subjclass[2020]{Primary 14J28, 14L24; Secondary 14J10}

\keywords{Polarized K3 surfaces, Chow stability,
geometric invariant theory, ADE singularities,
moduli spaces}

\begin{document}
\begin{abstract}
The asymptotic Chow stability of smooth polarized K3 surfaces
is known through the existence of Ricci-flat metrics and
Donaldson's theorem. We give a purely algebraic proof and
extend the result to polarized K3 surfaces with at worst ADE
singularities.
The proof gives a uniform effective bound for the exponent,
depending only on the degree of the polarization.
Following the approach of Mumford and Morrison, we establish
quantitative estimates for linear systems and their contributions
to Chow weights, without restrictions on the Picard number.
\end{abstract}
\maketitle

\section{Introduction}
\label{sec:introduction}

The study of GIT stability provides a natural approach to constructing
moduli spaces and their compactifications. For curves, the classical
results of Mumford and Gieseker establish the stability of sufficiently
high pluricanonical embeddings and yield a GIT construction of the
moduli space of stable curves \cite{Mumford,Gieseker}.
For surfaces, the picture is more varied.
Gieseker established asymptotic stability for canonically polarized
surfaces \cite{GiesekerSurfaces}, while Chow stability for
anticanonically polarized surfaces exhibits more subtle behavior (cf.~\cite{LeeLiSturmWang}).

The Ricci-flat case provides another natural setting for these
questions. For abelian surfaces, Chow stability of embeddings
defined by complete very ample linear systems admits a classical
algebraic proof using theta groups and Kempf's results \cite{Kempf}.
For projective K3 surfaces, the situation is more involved:
\begin{itemize}
\item In \cite{ViehwegIII,Viehweg}, Viehweg constructed quasi-projective moduli spaces in arbitrary
degree, using weak positivity and suitable GIT linearizations.

\item For polarized K3 surfaces of small degree, the following results are known:
 \begin{itemize}
     \item Shah gave GIT constructions in degrees $2$ and $4$
\cite{Shah2,Shah4}.
\item In \cite{LiTian}, Li and Tian proved GIT stability for complete-intersection K3
surfaces with at worst ADE singularities in degrees $6$ and $8$.
 \end{itemize}

\item In \cite{MorrisonK3}, Morrison proved asymptotic Hilbert stability for K3 surfaces
of Picard number one embedded by primitive polarizations of degree at
least $12$.

\item In \cite{Donaldson}, Donaldson used Ricci-flat K\"ahler metrics and the asymptotic
expansion of Bergman kernels to construct balanced embeddings of
smooth polarized K3 surfaces for sufficiently high powers of the
polarization. 
\end{itemize}

As mentioned in \cite{HuybrechtsK3}, Donaldson's approach does not identify the
ample line bundles on the moduli space arising from the construction,
nor does it clarify their relation to the ample line bundles obtained
from other constructions. A purely algebraic proof of asymptotic Chow stability remains
a natural problem. More generally,
R.~Pandharipande asked whether asymptotic Chow stability also holds
for polarized K3 surfaces with ADE singularities.

In this paper, we solve this problem via a purely algebraic method. More precisely, let $X$ be a complex projective K3 surface with at worst ADE
singularities, and let $L$ be an ample line bundle with $L^2=2d$.
For \(r\geq3\), the complete linear system \(|rL|\) defines an
embedding
\[
  X\hookrightarrow\PP^{N-1},
\]
where \(N=h^0(X,rL)=dr^2+2\).
\begin{theorem}\label{thm:main}
For every positive integer \(d\), there exists an integer
\(r_0=r_0(d)\) such that, for every complex polarized K3 surface
\((X,L)\) of degree \(2d\) with at worst ADE singularities,
the embedding defined by \(|L^{\otimes r}|\) is Chow stable
for all \(r\ge r_0\).
\end{theorem}

Our proof gives the explicit choice
\[
 r_0(d)=2^{31}(2d+1)^3.
\]
For smooth K3 surfaces, one may instead take
\(r_0(d)=2^{24}(2d+1)^3\).
These bounds are not expected to be optimal.
   
\begin{remark}
   Chow stability also implies asymptotic Hilbert stability
\cite{MorrisonGIT}. For each fixed $(X,L)$ and $r\geq r_0(d)$, the $q$-th Hilbert point of the
embedding defined by $L^{\otimes r}$ is stable for all sufficiently
large $q$. 
\end{remark}

Let $\cF_{2d}$ be the moduli stack of primitively polarized K3
surfaces of degree $2d$ with at worst ADE singularities, and let
$F_{2d}$ be its coarse moduli space. By the global Torelli theorem,
$F_{2d}$ is an arithmetic quotient of a type IV Hermitian symmetric
domain, hence is quasi-projective by the Baily--Borel theorem. By the standard construction of moduli spaces from embedded
polarized K3 surfaces, applied to the ADE locus
(cf.~\cite[Chapter~5, \S1.4 and \S2]{HuybrechtsK3}),
Theorem~\ref{thm:main} gives the following consequence.

\begin{corollary}
\label{cor:chow-compactification}
Fix an integer $r\ge r_0(d)$, and let $Z_{d,r}$ be the
reduced closure in the Chow variety of the locus of
primitively polarized K3 surfaces $(X,L)$ of degree $2d$
with at worst ADE singularities, embedded by $|L^{\otimes r}|$.
Then the Chow GIT quotient
\[
 Z_{d,r}\mathbin{/\mkern-6mu/}\operatorname{SL}_{dr^2+2}
\]
is a projective compactification of $F_{2d}$.
\end{corollary}

\subsection{Strategy of the Proof}
Our proof follows the algebraic approach of Morrison
\cite{MorrisonK3}, using Mumford's intersection-theoretic
estimates for Chow weights \cite{Mumford}. The main difficulty
in passing from Picard number one to arbitrary Picard number
is to control linear systems composed with a pencil and
those for which Clifford's inequality is nearly an equality.
We bound the total contribution of the former and show that
the latter, when their moving parts have positive
self-intersection, occur only \(O(r)\) times in any strict
flag of subspaces of \(H^0(X,rL)\). These estimates allow us
to prove Chow stability without imposing a condition on the
Picard number.

\subsection{Organization of the paper}

In Section~\ref{sec:main-argument}, we use Mumford's numerical
criterion to reduce Chow stability to an inequality comparing
the sum of the weights with an integral associated with the flag.
Section~\ref{sec:subsystem-estimates} establishes the required
weight estimates, treating moving parts of square zero and
positive square separately. In Section~\ref{sec:weighted-stability},
we combine these estimates to prove Theorem~\ref{thm:main}
and obtain an explicit bound depending only on the degree.

\subsection{Acknowledgements}
The author is supported by NSFC grant 12425105 and Shanghai Pilot Program for Basic Research (No. 21TQ001). He is also a member of LMNS. 

\subsection*{AI disclosure}
The author used an LLM to assist in developing the proof. It provided many useful suggestions concerning the function used in the integral
estimate in Section~2, several estimates in Section~3, and,
in particular, the subspace counting argument in Subsection~3.4.
The author takes full responsibility for the mathematical content
of the paper.

\section{Mumford's Numerical Criterion}
\label{sec:main-argument}

In this section, we recall Mumford's numerical criterion and reduce Chow stability
to an inequality comparing a sum of weights with an integral.

\subsection{Weighted Flags and the Degree of Contact}

Throughout, let \((X,L)\) be a complex polarized K3 surface with
at worst ADE singularities and \(L^2=2d\).
Fix \(r\ge3\), and write
\(N=h^0(X,rL)=dr^2+2\). The complete linear system \(|rL|\)
embeds \(X\) in \(\PP^{N-1}\), as recalled in
\cite{HuybrechtsK3} and \cite[\S5--6.1]{Morrison}.
By the Hilbert--Mumford criterion, Chow stability is tested by
one-parameter subgroups. We express their weights in terms of
Mumford's degree of contact.

\begin{definition}[Weighted flags and degree of contact]
\label{def:weighted-flag}
Let \(\lambda\) be a nontrivial one-parameter subgroup of
\(\operatorname{SL}(H^0(X,rL))\), with eigenbasis
\(z_0,\ldots,z_{N-1}\). After translating and rescaling its weights,
write
\[
 0=u_0\le\cdots\le u_{N-1}=1,\qquad
 V_j=\langle z_0,\ldots,z_j\rangle,\qquad
 W=\sum_{j=0}^{N-1}u_j.
\]
These data define a weighted flag \(\mathcal F\).
By its weight filtration, we mean the subspaces \(V_j\) at
the ends of the maximal blocks of equal weights. Let
\[
 \mathcal I_j=\operatorname{im}\bigl(
 V_j\otimes\mathcal O_X(-rL)\longrightarrow\mathcal O_X\bigr)
\]
be the associated base ideals. Choose a positive integer \(T\) with
\(Tu_j\in\mathbb Z\), and let \(p:X\times\mathbb A^1\to X\)
be the projection. Put
\[
 \mathfrak J=\sum_{j=0}^{N-1}t^{Tu_j}p^*\mathcal I_j,
 \qquad \mathcal L=p^*\mathcal O_X(rL).
\]
The normalized degree of contact \(e(\mathcal F)\) is defined, for
\(k\gg0\), by
\begin{equation}\label{eq:degree-contact-definition}
 \chi\!\left(
 \mathcal L^{\otimes k}/\mathfrak J^k\mathcal L^{\otimes k}
 \right)
 =\frac{T e(\mathcal F)}{6}\,k^3+O(k^2).
\end{equation}
Since \(t^T\in\mathfrak J\), the quotient has proper support.
The polynomial expansion follows from
\cite[Proposition~2.1 and Definition~2.2]{Mumford}.
Replacing \(T\) by \(qT\) amounts to the finite flat substitution
\(t\mapsto t^q\), which multiplies the Euler characteristic by \(q\).
Thus \(e(\mathcal F)\) is independent of \(T\).
\end{definition}

\begin{theorem}[Mumford {\cite[Theorem~2.9]{Mumford}}]
\label{thm:mumford-criterion}
The complete embedding by \(|rL|\) is Chow stable if and only if
every normalized weighted flag satisfies
\begin{equation*}
 e(\mathcal F)<\frac{3(rL)^2}{N}W.
\end{equation*}
\end{theorem}
We next bound the degree of contact using the moving parts
of the flag.

\subsection{The Degree of Contact and Moving Parts}
\label{subsec:moving-parts}

Let \(\pi:\widetilde X\to X\) be the minimal resolution.
For a nonzero subspace \(V\subseteq H^0(X,rL)\), choose a sequence
of point blowups \(\sigma:Y\to\widetilde X\) resolving its base ideal.
On \(Y\), write \(A=B+D\), where \(A\) represents the pullback
of \(rL\), and \(D\) is the fixed part of the pulled-back linear
system defined by \(V\). We call \(B\) the moving part associated
with \(V\). The induced evaluation map
\[
 V\otimes\mathcal O_Y\longrightarrow\mathcal O_Y(B)
\]
is surjective, so \(B\) is globally generated and nef.
For a flag \((V_j)\), we choose a common resolution and write
\(A=B_j+D_j\) for the corresponding decompositions.

\begin{proposition}\label{prop:contact-intersection-estimate}
Let \(\mathcal F\) be a normalized weighted flag as in
Definition~\ref{def:weighted-flag}. There exist a sequence of point
blowups \(\sigma:Y\to\widetilde X\) and divisors
\(A,B_0,\ldots,B_{N-1}\) on \(Y\) such that
\[
 \mathcal O_Y(A)=(\pi\circ\sigma)^*\mathcal O_X(rL),
\]
where \(B_j\) is the moving part associated with \(V_j\).
In particular,
the \(B_j\) are nef and satisfy
\begin{equation*}
 0=B_0\le B_1\le\cdots\le B_{N-1}=A,
 \qquad h^0(Y,B_j)\ge j+1.
\end{equation*}
Moreover, for every subsequence
\(0=i_0<\cdots<i_\ell=N-1\), one has
\begin{equation}\label{eq:mumford-global-path}
 e(\mathcal F)\le
 3A^2-\sum_{\nu=0}^{\ell-1}
 (u_{i_{\nu+1}}-u_{i_\nu})
 \bigl(B_{i_\nu}^2+B_{i_\nu}\cdot B_{i_{\nu+1}}
                       +B_{i_{\nu+1}}^2\bigr).
\end{equation}
\end{proposition}

\begin{proof}
Resolve the base ideals \(\mathcal I_j\) simultaneously by point
blowups \(\sigma:Y\to\widetilde X\), and write
\[
 \mathcal I_j\mathcal O_Y=\mathcal O_Y(-D_j).
\]
Let \(A\) be the divisor of the pullback of \(z_0\), and set
\(B_j=A-D_j\). The evaluation maps give surjections
\[
 V_j\otimes\mathcal O_Y\longrightarrow\mathcal O_Y(B_j).
\]
Thus the \(B_j\) are globally generated and nef, with
\(h^0(Y,B_j)\ge j+1\).
Since \(\mathcal I_j\subseteq\mathcal I_{j+1}\), we have
\(D_j\ge D_{j+1}\). Moreover, \(D_0=A\) and \(D_{N-1}=0\),
because \(V_0=\langle z_0\rangle\) and \(|rL|\) is
base-point-free. This proves the assertions about the \(B_j\).

Fix a subsequence \(0=i_0<\cdots<i_\ell=N-1\), and suppose
first that its weights are distinct.
Omitting subspaces from the flag can only increase the degree
of contact. Applying Mumford's estimates to this subsequence
\cite[Propositions~4.9--4.10 and Corollary~4.11]{Mumford},
we obtain
\[
 e(\mathcal F)\le
 \sum_{\nu=0}^{\ell-1}(u_{i_{\nu+1}}-u_{i_\nu})
 \bigl(e_A(\mathcal I_{i_\nu})
 +e_A(\mathcal I_{i_\nu},\mathcal I_{i_{\nu+1}})
 +e_A(\mathcal I_{i_{\nu+1}})\bigr),
\]
where \(e_A\) denotes the global multiplicities and mixed
multiplicities with respect to \(A\).
By \cite[Proposition~4.8 and Lemma~5.6]{Mumford},
\[
 e_A(\mathcal I_i)=A^2-B_i^2,\qquad
 e_A(\mathcal I_i,\mathcal I_j)=A^2-B_i\cdot B_j,
\]
with the mixed formula obtained by polarization.
Substituting these identities and using
\(\sum_{\nu=0}^{\ell-1}(u_{i_{\nu+1}}-u_{i_\nu})=1\)
gives \eqref{eq:mumford-global-path}.

For repeated weights, keep only the last subspace in each
block of equal weights in the chosen subsequence. Indeed, if
\(i<j\) and \(u_i=u_j\), then \(V_i\subseteq V_j\), so \(V_i\)
contributes no additional sections at that weight. Removing it
therefore leaves the degree of contact of the subsequence unchanged.
The remaining weights are distinct and still run from \(0\) to
\(1\), so the preceding calculation applies.

To compare the resulting bound with the original sum, note
that for \(i\le i'\) and \(j\le j'\),
\[
 B_{i'}\cdot B_{j'}-B_i\cdot B_j
 =(B_{i'}-B_i)\cdot B_{j'}+B_i\cdot(B_{j'}-B_j)\ge0,
\]
since the differences are effective and the \(B_j\) are nef.
Replacing each block by its last subspace therefore increases
the intersection terms being subtracted, while terms between
equal weights vanish. This proves
\eqref{eq:mumford-global-path} in general.
\end{proof}

\subsection{A Convex Function Associated with the Flag}
Fix \(Y\) and \(B_j\) as in
Proposition~\ref{prop:contact-intersection-estimate}, and put
\[
 a=\frac{A^2}{2}=dr^2=N-2,\qquad b_j=B_j^2.
\]
We associate with the weighted flag a function that records the
largest self-intersection attainable by convex combinations of
the moving parts, subject to a bound on their average weight.
Convex combinations are needed to control the moving parts
of square zero. Although such a divisor has zero
self-intersection, its mixed intersection with a nef divisor
of positive square can contribute to the square of their
convex combination. We use this observation to estimate the
contribution of pencils in Subsection~\ref{sec:contact-pencil}.

\begin{definition}
\label{def:volume-envelope}
For \(0\le\rho\le1\), set
\begin{equation*}
 \phi(\rho)=
 \max\left\{
 \sqrt{\frac{Q^2}{2a}}:\quad
 Q=\sum_j\lambda_jB_j,\quad \lambda_j\ge0,\quad
 \sum_j\lambda_j=1,\quad
 \sum_j\lambda_ju_j\le\rho
 \right\},
\end{equation*}
and define
\[
 h(y)=\min\{\rho\in[0,1]:\phi(\rho)\ge y\},
 \qquad 0\le y\le1.
\]
\end{definition}

\begin{proposition}
\label{prop:envelope-properties}
The function \(\phi\) is continuous, nondecreasing and concave,
and satisfies
\[
 \rho\le\phi(\rho)\le1,\qquad \phi(1)=1.
\]
Its generalized inverse \(h\) is continuous, nondecreasing and
convex. Moreover,
\begin{equation}\label{eq:coordinate-envelope}
 h(0)=0,\qquad h(y)\le y,\qquad
 u_j\ge h\!\left(\sqrt{\frac{b_j}{2a}}\right).
\end{equation}
\end{proposition}

\begin{proof}
The admissible coefficient sets in the definition of \(\phi\)
are compact and increase with \(\rho\). Thus \(\phi\) is
nondecreasing. Concavity follows from
\[
 \sqrt{(P+Q)^2}\ge\sqrt{P^2}+\sqrt{Q^2}
\]
for nef divisor classes \(P,Q\), by the Hodge index theorem.

Since \(B_0=0\) and \(B_{N-1}=A\), the divisor \(\rho A\)
is admissible for \(\phi(\rho)\). On the other hand, every
admissible \(Q\) is nef and \(A-Q\) is effective, so
\(Q^2\le A^2=2a\). Hence
\[
 \rho\le\phi(\rho)\le1.
\]
Compactness gives continuity at \(0\), while these bounds give
continuity at \(1\). Concavity gives continuity in the interior.
The asserted properties of \(h\) now follow from the
corresponding properties of the generalized inverse of a
continuous, nondecreasing concave function.

Finally, \(B_j\) is admissible at \(\rho=u_j\), so
\[
 \phi(u_j)\ge\sqrt{\frac{b_j}{2a}}.
\]
The definition of \(h\) gives the last inequality in
\eqref{eq:coordinate-envelope}, while \(\phi(\rho)\ge\rho\)
gives \(h(y)\le y\).
\end{proof}

\begin{proposition}
\label{prop:contact-envelope}
For every weighted flag as above,
\begin{equation}
 e(\mathcal F)
 \le6a\int_0^1\bigl(1-\phi(\rho)^2\bigr)\,d\rho
 =6a\int_0^1h(\sqrt{x})\,dx.
 \label{eq:contact-envelope-bound}
\end{equation}
\end{proposition}
\begin{proof}
Choose a rational partition
\(0=\rho_0<\cdots<\rho_\ell=1\).
For each \(i\), choose rational coefficients \(\lambda_{ij}\)
such that
\[
 Q_i=\sum_j\lambda_{ij}B_j,\qquad
 \lambda_{ij}\ge0,\qquad \sum_j\lambda_{ij}=1,\qquad
 \sum_j\lambda_{ij}u_j\le\rho_i,
\]
and \(Q_i^2/(2a)\) approximates \(\phi(\rho_i)^2\).
Such approximations exist because the admissible coefficient
sets are rational polytopes. At the endpoint, choose \(Q_\ell=A\).
Take \(q\) sufficiently divisible that all \(q\lambda_{ij}\)
and \(q\rho_i\) are integral. The coefficient ideals of
\(\mathfrak J^q\) are
\[
 \mathcal K_{q,\nu}
 =\sum_{\substack{\sum_j n_j=q\\T\sum_jn_ju_j\le\nu}}
       \prod_j\mathcal I_j^{n_j}.
\]
Let \(W_i\subseteq H^0(X,qrL)\) be generated by products of
\(q\) basis sections of total weight at most \(q\rho_i\).
Its base ideal is \(\mathcal K_{q,qT\rho_i}\).

Resolve these ideals on a common resolution dominating \(Y\),
and denote the moving parts by \(C_i\). We use the same notation
for the pullbacks of \(A\), \(B_j\), and \(Q_i\).
Since \(V_j\) generates \(\mathcal O_Y(B_j)\), products with
\(q\lambda_{ij}\) factors from each \(V_j\) generate the line
bundle associated with \(qQ_i\). They form a subsystem of
\(W_i\), so \(C_i-qQ_i\) is effective. Since these divisors
are nef,
\[
 C_i^2\ge q^2Q_i^2,\qquad
 C_i\cdot C_{i+1}\ge q^2Q_i\cdot Q_{i+1}.
\]
At \(\rho_\ell=1\), all products are allowed and generate
\(\mathcal O_X(qrL)\), hence \(C_\ell=qA\).

The sheaves \(\mathcal K_{q,\nu}\otimes\mathcal O_X(qrL)\)
and \(\mathfrak J^q\mathcal L^{\otimes q}\) are generated by
the corresponding products. Mumford's surface estimate
\cite[Proposition~4.10 and Corollary~4.11]{Mumford} therefore applies.
For fixed \(q\), equation~\eqref{eq:degree-contact-definition} gives
\[
 \chi\!\left(
 \mathcal L^{\otimes qk}/\mathfrak J^{qk}\mathcal L^{\otimes qk}
 \right)
 =\frac{Tq^3e(\mathcal F)}{6}\,k^3+O(k^2).
\]
Thus the left-hand side of Mumford's estimate is
\(Tq^3e(\mathcal F)\), while the mixed multiplicities are
\[
 \begin{aligned}
 e_{qA}(\mathcal K_{q,qT\rho_i})&=q^2A^2-C_i^2,\\
 e_{qA}(\mathcal K_{q,qT\rho_i},\mathcal K_{q,qT\rho_{i+1}})
 &=q^2A^2-C_i\cdot C_{i+1}.
 \end{aligned}
\]
Dividing the resulting estimate by \(Tq^3\) gives
\[
 e(\mathcal F)\le6a-
 \sum_{i=0}^{\ell-1}(\rho_{i+1}-\rho_i)
 \bigl(Q_i^2+Q_i\cdot Q_{i+1}+Q_{i+1}^2\bigr).
\]
By Hodge index, \(Q_i\cdot Q_{i+1}\ge\sqrt{Q_i^2Q_{i+1}^2}\).
We now let the approximation errors tend to zero and refine the
partition. By continuity of \(\phi\), we obtain
\[
 e(\mathcal F)\le6a\int_0^1(1-\phi(\rho)^2)\,d\rho.
\]
Finally,
\(h(\sqrt{x})\le\rho\) if and only if \(x\le\phi(\rho)^2\).
Fubini's theorem gives the asserted equality of integrals.
\end{proof}

Combining Proposition~\ref{prop:contact-envelope} with
Mumford's criterion, we obtain the following sufficient condition.

\begin{proposition}\label{prop:scalar-chow-criterion}
The embedding defined by \(|rL|\) is Chow stable if every
normalized weighted flag \(\mathcal F\) satisfies
\begin{equation}\label{eq:scalar-sufficient}
 N\int_0^1h(\sqrt{x})\,dx<W,
\end{equation}
where \(h\) and \(W\) are associated with \(\mathcal F\).
\end{proposition}

Section~\ref{sec:subsystem-estimates} gives geometric estimates
for the weights, which are used in Section~\ref{sec:weighted-stability}
to verify \eqref{eq:scalar-sufficient}.

\section{Weight Estimates for Linear Subsystems}
\label{sec:subsystem-estimates}

To verify \eqref{eq:scalar-sufficient}, we compare values of \(h\)
with the weights \(u_j\). We treat moving parts of square zero
and positive square separately, then collect the comparisons in
Proposition~\ref{prop:weight-estimates}. Throughout, \(B\) denotes
the moving part associated with a subspace \(V\), as defined in
Subsection~\ref{subsec:moving-parts}.

\subsection{Linear Systems Composed with a Pencil}
\label{sec:contact-pencil}

For the moving parts with \(B^2=0\), we bound the total error
in the comparison with the weights.
Let \(P\) be the largest \(j\in\{1,\ldots,a\}\) with \(B_j^2=0\),
setting \(P=0\) if there is none. We write \((x)_+=\max\{x,0\}\).

\begin{lemma}
\label{lem:localized-pencil-correction}
Put \(s=16(2d+1)\). Assume
\begin{equation}
 r\ge32(s+1),\qquad a\ge16s,\qquad
 \frac{s(1+\log(2dr))}{r}\le1.
 \label{eq:pencil-size-assumptions}
\end{equation}
Then
\begin{equation}\label{eq:localized-pencil-cost}
\begin{aligned}
 \sum_{j=1}^{P}
 \left(h\!\left(\sqrt{\frac{j+s}{a}}\right)-u_j\right)_+
 &\le16\left[2d+\frac{s(1+\log(2dr))}{r}\right]h(1/2)\\
 &\le s\,h(1/2).
\end{aligned}
\end{equation}
\end{lemma}

\begin{proof}
Suppose \(P\ge1\), since otherwise there is nothing to prove.
The inequalities \(0\le B_i\le B_j\) and nefness imply that
\(B_j^2\) is nondecreasing. For \(1\le i<j\le P\),
\[
 0\le B_i\cdot B_j
 =B_j^2-B_j\cdot(B_j-B_i)\le0.
\]
Thus the nonzero isotropic divisors are proportional and define
the same fibration \(f:Y\to\mathbb P^1\), since \(q(Y)=0\).
For its general connected fiber \(F\), we have
\[
 B_j\sim q_jF,\qquad h^0(Y,B_j)=q_j+1,\qquad q_j\ge j.
\]
Put \(\ell=A\cdot F\). Since \(A\) is nef and big and is
\(r\) times an integral divisor, \(\ell\ge r\).
The effectivity of \(A-q_PF\) gives
\begin{equation}\label{eq:pencil-prefix-size}
 P\le q_P\le2a/\ell\le2dr.
\end{equation}
In particular, the assumptions imply
\begin{equation}\label{eq:pencil-low-grid}
 P+s\le2a/r+s\le a/8+a/16<a/4.
\end{equation}

Set \(\theta=h(1/2)\). Choose an admissible convex combination
\(Q\) at weight \(\theta\) with \(Q^2\ge a/2\).
If necessary, replace \(Q\) by \(\lambda Q+(1-\lambda)B_0\)
for a suitable \(0<\lambda\le1\). This preserves admissibility
and gives \(Q^2=a/2\). Since \(Q\) is a convex combination
of the \(B_j\), the divisor \(A-Q\) is effective.
For a nonzero nef divisor \(G\) with \(G^2=0\), Hodge index
on the orthogonal complement of the plane spanned by \(A,G\)
gives, for every real divisor class \(B'\),
\begin{equation}\label{eq:isotropic-hodge-bound}
 (B')^2\le
 \frac{2(A\cdot B')(B'\cdot G)}{A\cdot G}
 -\frac{2a(B'\cdot G)^2}{(A\cdot G)^2}.
\end{equation}
Applying this to \(Q,F\), and using \(A\cdot Q\le2a\), yields
\[
 a/2=Q^2\le\frac{4a(F\cdot Q)}{\ell},
 \qquad F\cdot Q\ge\ell/8.
\]

For \(1\le j\le P\), put
\[
 t_j=\frac{16(j+s)}{q_j\ell}
 \le\frac{16(s+1)}r\le\frac12.
\]
The convex combination \((1-t_j)B_j+t_jQ\) is admissible
at weight \((1-t_j)u_j+t_j\theta\), and
\[
 \bigl((1-t_j)B_j+t_jQ\bigr)^2
 \ge2t_j(1-t_j)q_j(F\cdot Q)
 \ge\frac{t_jq_j\ell}{8}=2(j+s).
\]
Consequently,
\[
 \left(h\!\left(\sqrt{\frac{j+s}{a}}\right)-u_j\right)_+
 \le t_j(\theta-u_j)_+\le t_j\theta.
\]
Finally, \(q_j\ge j\), \eqref{eq:pencil-prefix-size}, and the
harmonic sum bound give
\[
 \begin{aligned}
 \sum_{j=1}^P t_j
 &\le\frac{16}{\ell}\bigl(P+s(1+\log P)\bigr)\\
 &\le16\left[2d+\frac{s(1+\log(2dr))}{r}\right]\le s.
 \end{aligned}
\]
Multiplying by \(\theta\) proves the assertion.
\end{proof}

\subsection{Linear Systems with Two-Dimensional Image}
\label{subsec:positive-square}

When \(B^2>0\), Clifford's inequality controls \(\dim V\)
in terms of \(B^2\). Combined with \eqref{eq:coordinate-envelope},
this gives a pointwise comparison between \(h\) and the weights.
For a globally generated divisor \(B\) on \(Y\) and a subspace
\(V\subseteq H^0(Y,B)\) generating \(B\), put
\begin{equation*}
 \kappa(B,V)=B^2-2(\dim V-2),\qquad
 \kappa(B)=B^2-2(h^0(Y,B)-2).
\end{equation*}

\begin{proposition}\label{prop:clifford-comparison}
Suppose \(B^2>0\). Then
\begin{equation}\label{eq:moving-clifford-bound}
 h^0(Y,B)\le\frac{B^2}{2}+2,
 \qquad \kappa(B,V)\ge\kappa(B)\ge0.
\end{equation}
If \(V\) belongs to the filtration associated with \(\lambda\),
\(B\) is its moving part, and
\(\kappa(B,V)\ge2s+2\) for an integer \(s\ge1\), then every
index \(j\) in the corresponding block of equal weights satisfies
\begin{equation}\label{eq:clifford-shifted-comparison}
 u_j\ge h\!\left(\sqrt{\frac{\min(j+s,a)}{a}}\right).
\end{equation}
\end{proposition}

\begin{proof}
Let \(C\in|B|\) be a general member. Since \(B^2>0\) and
\(B\) is globally generated, \(C\) is smooth and irreducible.
Moreover, \(q(Y)=0\), and the restriction sequence gives
\[
 h^0(C,\mathcal O_C(B))=h^0(Y,B)-1.
\]
By adjunction,
\[
 K_C\otimes\mathcal O_C(-B)=\mathcal O_C(K_Y).
\]
Since \(K_Y\) is effective and \(C\) is general,
\(\mathcal O_C(B)\) is special. Clifford's inequality gives
\[
 h^0(C,\mathcal O_C(B))\le\frac{B^2}{2}+1.
\]
The remaining inequalities follow from
\[
 \kappa(B,V)=\kappa(B)+2\bigl(h^0(Y,B)-\dim V\bigr).
\]

For the weight comparison, put \(v=\dim V\). The hypothesis
\(\kappa(B,V)\ge2s+2\) gives
\[
 B^2\ge2(v+s-1).
\]
Since \(V\) is the last subspace in its block of equal weights,
\(u_j=u_{v-1}\) for every index \(j\) in that block, and
\(j\le v-1\). Applying \eqref{eq:coordinate-envelope} and
using the monotonicity of \(h\), we obtain
\[
 u_j=u_{v-1}
 \ge h\!\left(\sqrt{\frac{B^2}{2a}}\right)
 \ge h\!\left(\sqrt{\frac{\min(j+s,a)}{a}}\right).
\]
\end{proof}

\subsection{Bounds for Base-Point Multiplicities}
\label{subsec:base-point-multiplicities}

The shifted comparison above applies when \(\kappa(B,V)\) is
large. To count the remaining subspaces, we first bound the
base-point multiplicities when \(\kappa(B)\) is small.

\begin{lemma}\label{lem:base-point-multiplicity-bound}
Let \(\eta\ge0\) be an integer, and let \(B\) be a globally
generated divisor on \(Y\) such that \(A-B\) is effective.
If
\[
 \kappa(B)\le\eta,\qquad
 B^2>\max\{8dr,\,4\eta+8\},
\]
then \(K_Y\cdot B\le3\eta\).
\end{lemma}

\begin{proof}
Put \(b=B^2\), \(m=h^0(Y,B)-2\), \(\kappa=\kappa(B)\), and
\(\delta=K_Y\cdot B\). The assumptions imply \(m\ge\eta+3\).
The morphism defined by \(|B|\) has nondegenerate image
\(T\subset\PP^{m+1}\). If its degree is \(e\), then
\[
 em\le e\deg T=b=2m+\kappa<3m,
\]
so \(e=1\) or \(2\).

If \(e=1\), the complete linear series on a general member
\(C\in|B|\) has degree \(2m+\kappa\) and dimension \(m\), since
\(q(Y)=0\). Castelnuovo's bound \cite[\S1]{Harris} gives
\begin{equation}\label{eq:genus-complete-moving-curve}
 m+1+\frac{\kappa+\delta}{2}=g(C)\le m+2\kappa+1.
\end{equation}
Hence \(\delta\le3\kappa\le3\eta\).

Suppose \(e=2\). Let \(Q\) be the minimal resolution of the
normalization of \(T\), and let \(D\) be the pullback of its
hyperplane divisor. Resolving the rational map gives morphisms
\[
 p:Y'\longrightarrow Y,\qquad f:Y'\longrightarrow Q,
 \qquad p^*B=f^*D,
\]
with \(p\) birational and \(\deg f=2\). Pullback of sections and
of one-forms gives
\[
 h^0(Q,D)=m+2,\qquad D^2=m+\kappa/2,\qquad q(Q)=0.
\]
Thus \(|D|\) is complete and birational onto \(T\).
Castelnuovo's bound gives \(g\le\kappa/2\) for the genus of a
general member of \(|D|\). In particular,
\(K_Q\cdot D=2g-2-D^2<0\), so \(Q\) is rational by
Castelnuovo's rationality criterion.

Apply \cite[Propositions~4.1 and 5.1]{CDML}.
If \(g>0\), then \(D^2\le4g+4\le2\eta+4\), unless \(T\)
is the Veronese surface of degree \(9\). Both alternatives
contradict \(b>\max\{8dr,4\eta+8\}\), since \(8dr\ge24\).
If \(g=0\), then \(T\) is a scroll, since \(D^2>4\).
Its ruling on \(Q\) has general fiber \(F\) with
\(F^2=0\) and \(D\cdot F=1\).

Set \(G=f^*F\). Then \(G\) is nonzero and nef, with
\(G^2=0\) and \(p^*B\cdot G=2\).
Continue to denote the pullback of \(A\) by \(A\).
Since \(A\) is nef and big and divisible by \(r\),
\(A\cdot G\ge r\). Moreover, \(A\cdot p^*B\le A^2=2a\).
Applying \eqref{eq:isotropic-hodge-bound}, we obtain
\[
 b\le\frac{2(A\cdot p^*B)(p^*B\cdot G)}{A\cdot G}
 \le\frac{8a}{r}=8dr,
\]
a contradiction.
\end{proof}

\subsection{Counting Subspaces in a Flag}
\label{subsec:counting-subspaces}

We now show that only \(O(r)\) subspaces in a strict flag can
have positive square and bounded \(\kappa(B,V)\).

Fix an integral ample divisor \(H\) on \(\widetilde X\), and put
\(c=H\cdot\pi^*L\). For a strict flag
\[
 0\ne V_0\subsetneq\cdots\subsetneq V_n\subseteq H^0(X,rL),
\]
let \(B_i\) be the moving parts on a common resolution
\(\sigma:Y\to\widetilde X\) of the base ideals.

\begin{proposition}\label{prop:uniform-subspace-count}
For every integer \(\eta\ge0\), there is a constant \(C>0\),
depending only on \(d,c,\eta\), such that
\[
 \#\{i:B_i^2>0,\ \kappa(B_i,V_i)\le\eta\}\le Cr.
\]
The bound is independent of the flag and its resolution.
\end{proposition}

\begin{proof}
Put
\[
 R=\max\{16dr,\,4\eta+8\}.
\]
For \(0<B_i^2\le R\), Clifford's inequality gives
\[
 \dim V_i\le h^0(Y,B_i)\le R/2+2.
\]
Since the dimensions in a strict flag are distinct, there
are at most \(\lfloor R/2\rfloor+2\) such indices.

Let \(F_i\) be the fixed divisor of \(V_i\) on
\(\widetilde X\), and set \(M_i=r\pi^*L-F_i\).
The inclusions \(V_i\subseteq V_j\) imply
\(F_i\ge F_j\), hence \(M_i\le M_j\), for \(i<j\).
Moreover,
\[
 0\le H\cdot M_i\le rc.
\]
Every strict increase of \(M_i\) increases its integral
\(H\)-degree, so at most \(rc+1\) distinct divisors \(M_i\)
occur.

Now consider an index with \(B_i^2>R\) and
\(\kappa(B_i,V_i)\le\eta\).
Writing \(E_\alpha\) for the total transforms of the
exceptional curves, we have
\[
 B_i=\sigma^*M_i-\sum_\alpha m_{i,\alpha}E_\alpha,
 \qquad m_{i,\alpha}\ge0.
\]
The blowup formulas and
Lemma~\ref{lem:base-point-multiplicity-bound} give
\[
 \sum_\alpha m_{i,\alpha}=K_Y\cdot B_i\le3\eta,
 \qquad
 M_i^2-B_i^2=\sum_\alpha m_{i,\alpha}^2\le9\eta^2.
\]
Since
\[
 \dim V_i=\frac{B_i^2-\kappa(B_i,V_i)}2+2,
 \qquad 0\le\kappa(B_i,V_i)\le\eta,
\]
all such indices with a fixed \(M_i=M\) satisfy
\[
 \frac{M^2-9\eta^2-\eta}{2}+2
 \le\dim V_i\le\frac{M^2}{2}+2.
\]
There are therefore at most
\((9\eta^2+\eta)/2+1\) such indices for each \(M\).
Together with the preceding counts, this gives
\begin{equation}\label{eq:explicit-subspace-count}
 \#\{i:B_i^2>0,\ \kappa(B_i,V_i)\le\eta\}
 \le
 \left\lfloor\frac R2\right\rfloor+2+
 \left(\frac{9\eta^2+\eta}{2}+1\right)(rc+1).
\end{equation}
For fixed \(d,c,\eta\), the right-hand side is bounded by \(Cr\).
All estimates are independent of the chosen resolution.
\end{proof}

\subsection{Comparison with the Weights}
\label{subsec:positive-square-weights}

Put \(s=16(2d+1)\), and write
\[
 x_j=\sqrt{j/a}\quad(0\le j\le a),\qquad
 y_j=x_{\min(j+s,a)}\quad(1\le j\le a).
\]
As in Subsection~\ref{sec:contact-pencil}, \(P\) is the largest
index in \(\{1,\ldots,a\}\) for which \(B_P^2=0\), or \(0\)
if there is no such index. 

\begin{proposition}\label{prop:weight-estimates}
Assume \eqref{eq:pencil-size-assumptions}. For the indices with
\(B_j^2=0\), the total error satisfies
\[
 \sum_{j=1}^{P}\bigl(h(y_j)-u_j\bigr)_+\le s\,h(1/2).
\]
For the remaining indices, there are sets
\(E_1\subseteq E_s\subseteq\{P+1,\ldots,a\}\) such that
\begin{equation}\label{eq:weight-lower-bounds}
 \begin{aligned}
 h(x_{j-1})&\le u_j &&(P<j\le a),\\
 h(x_j)&\le u_j &&(P<j\le a,\ j\notin E_1),\\
 h(y_j)&\le u_j &&(P<j\le a,\ j\notin E_s).
 \end{aligned}
\end{equation}
The exceptional sets satisfy
\begin{align}
 |E_1|&\le k,\label{eq:ordinary-exception-count}\\
 |E_s|&\le(s+1)k,\label{eq:shifted-exception-count}
\end{align}
where \(k\) is the number of subspaces \(V\) in the weight
filtration for which \(B^2>0\) and \(\kappa(B,V)\le2s+1\).
Moreover, \(k\le Cr\), with \(C\) depending only on \(d,c\).
\end{proposition}

\begin{proof}
The bound on \(k\) follows from
Proposition~\ref{prop:uniform-subspace-count}.
Lemma~\ref{lem:localized-pencil-correction} gives the
asserted estimate for \(j\le P\), since
\eqref{eq:pencil-low-grid} implies \(P+s<a\).

Consider a maximal block \(I\) of equal weights, and let
\(v-1\) be its last index. For \(j\in I\) with \(j>P\),
the divisor \(B_{v-1}\) has positive square. Put
\[
 b=B_{v-1}^2,\qquad
 \kappa=\kappa(B_{v-1},V_{v-1})=b-2(v-2).
\]
Clifford's inequality and \eqref{eq:coordinate-envelope}
give
\[
 u_j=u_{v-1}
 \ge h\!\left(\sqrt{\frac{b}{2a}}\right)
 \ge h(x_{j-1}),
\]
because \(b\ge2(v-2)\) and \(j\le v-1\).

If \(\kappa\ge2s+2\), Proposition~\ref{prop:clifford-comparison}
gives \(u_j\ge h(y_j)\) throughout the block.
For the remaining blocks, the same bound \(b\ge2(v-2)\)
gives
\[
 \begin{aligned}
 u_j&\ge h(x_j) &&\text{when }j\le v-2,\\
 u_j&\ge h(y_j) &&\text{when }j+s\le v-2.
 \end{aligned}
\]
Thus the first comparison can fail only at the last index
of the block, and the second only at its last \(s+1\) indices.

Define \(E_1\) and \(E_s\) to consist of these respective
indices in \(\{P+1,\ldots,a\}\), over the blocks with
\(b>0\) and \(\kappa\le2s+1\).
There are at most \(k\) such blocks, so
\eqref{eq:ordinary-exception-count} and
\eqref{eq:shifted-exception-count} follow.
The preceding comparisons prove
\eqref{eq:weight-lower-bounds}.
\end{proof}

\section{Proof of the Main Theorem}
\label{sec:weighted-stability}
\label{sec:numerical-comparison}
\label{sec:scalar-flags}

We compare the integral in \eqref{eq:scalar-sufficient} with
finite sums of values of \(h\), then use
Proposition~\ref{prop:weight-estimates} to bound those sums by
the weights. This proves the required inequality once the
exceptional sets are sufficiently small.

\subsection{Riemann Sum Estimates}

We first estimate the difference between the discrete sum and
the integral for the functions \((t-c)_+\). The piecewise linear interpolation reduces
the estimate for \(h\) to nonnegative linear combinations of these functions.

\begin{lemma}\label{lem:hinge-quadrature}
Let \(a\ge8\) be an integer, and put \(x_j=\sqrt{j/a}\).
For \(0\le c<1\), set \(\psi_c(t)=(t-c)_+=\max\{t-c,0\}\) and
\[
 Q_0[\psi_c]=\psi_c(1)+\sum_{j=1}^a \psi_c(x_j)
              -(a+2)\int_0^1\psi_c(\sqrt{x})\,dx.
\]
Then
\begin{align}
 Q_0[\psi_c]&>\frac1{12} &&(0\le c\le4/5),
 \label{eq:hinge-low-base}\\
 Q_0[\psi_c]&>1-c &&(4/5\le c<1).
 \label{eq:hinge-high-base}
\end{align}
\end{lemma}

\begin{proof}
Direct integration gives
\[
 \int_0^1\psi_c(\sqrt{x})\,dx=\frac23-c+\frac{c^3}{3}.
\]
Strict concavity of the square root gives the midpoint estimate
\[
 \begin{aligned}
 \sum_{j=1}^a\sqrt{j/a}
 &>\frac1{\sqrt a}\int_{1/2}^{a+1/2}\sqrt{t}\,dt\\
 &=\frac{2}{3\sqrt a}
   \bigl((a+1/2)^{3/2}-(1/2)^{3/2}\bigr)
 >\frac{2a}{3}+\frac12-\frac1{3\sqrt{2a}}.
 \end{aligned}
\]
For the convex function \(x\mapsto(c-\sqrt{x})_+\), the
trapezoidal estimate gives
\[
 \sum_{j=1}^a(c-x_j)_+
 -a\int_0^1(c-\sqrt{x})_+\,dx\ge-\frac c2.
\]
Using \(\psi_c(\sqrt{x})=\sqrt{x}-c+(c-\sqrt{x})_+\)
and adding these estimates, we obtain
\begin{equation}\label{eq:hinge-Q0}
 \sum_{j=1}^a\psi_c(x_j)
 -a\int_0^1\psi_c(\sqrt{x})\,dx
 >\frac12-\frac1{3\sqrt{2a}}-\frac c2
 \ge\frac5{12}-\frac c2.
\end{equation}
Consequently,
\[
 Q_0[\psi_c]>\frac1{12}+\frac c2-\frac{2c^3}{3}
 \ge\frac1{12}\qquad(0\le c\le4/5).
\]
For \(4/5\le c<1\), integration on the grid intervals yields
\[
 \sum_{j=1}^a\psi_c(x_j)
 -a\int_0^1\psi_c(\sqrt{x})\,dx
 =\int_{c^2}^1\frac{\{ax\}}{2\sqrt{x}}\,dx
 \ge\frac{1-c^2}{4}.
\]
Here \(\{ax\}\) is the fractional part of \(ax\).
Its average is \(1/2\) on each full grid interval and at least
\(1/2\) on the initial partial interval. Hence
\[
 Q_0[\psi_c]-(1-c)
 \ge\frac{1-c}{12}(8c^2+11c-13)>0.
\]
\end{proof}

\subsection{Summation of the Weight Estimates}
We now combine the weight estimates with the Riemann sum estimates
of the preceding subsection. The following proposition gives the
inequality required by Proposition~\ref{prop:scalar-chow-criterion}.

\begin{proposition}\label{prop:weighted-chow-inequality}
Let \(\mathcal F\) be a normalized weighted flag, and let
\(h\) be its associated function. Put \(s=16(2d+1)\), and
let \(k\) count the subspaces in its weight filtration with
\(B^2>0\) and \(\kappa(B,V)\le2s+1\).
If \eqref{eq:pencil-size-assumptions} holds and
\begin{equation}\label{eq:scalar-numerical-conditions}
 \frac{k}{a}\le
 \frac{s^2}{64(s+1)^3},
\end{equation}
then
\begin{equation}\label{eq:scalar-target}
 W>N\int_0^1h(\sqrt{x})\,dx.
\end{equation}
\end{proposition}
\begin{proof}
If \(h=0\), then \(W\ge u_{a+1}=1\) proves the assertion.
Assume that \(h\ne0\), and recall
\[
 x_j=\sqrt{j/a},\qquad y_j=x_{\min(j+s,a)}.
\]

We first reduce to a finite nonnegative linear combination of
the functions \(\psi_c\) from Lemma~\ref{lem:hinge-quadrature}.
Let \(\widetilde h\) be the piecewise linear interpolation of
\(h\) at \(x_0,\ldots,x_a\) and \(1/2\).
It is convex and nondecreasing, and convexity of \(h\) gives
\(h\le\widetilde h\). Moreover,
\[
 \widetilde h(x_j)=h(x_j),\qquad
 \widetilde h(y_j)=h(y_j),\qquad
 \widetilde h(1/2)=h(1/2),
\]
because each \(y_j\) is one of the points \(x_0,\ldots,x_a\).
Thus replacing \(h\) by \(\widetilde h\) leaves both the
pointwise weight bounds and the total error bound in
Proposition~\ref{prop:weight-estimates} unchanged.
It suffices to prove
\[
 W>N\int_0^1\widetilde h(\sqrt{x})\,dx.
\]

Since \(\widetilde h(0)=0\) and its slopes are nonnegative
and nondecreasing, we can write
\[
 \widetilde h=\sum_{\ell=1}^m\beta_\ell\psi_{c_\ell},
 \qquad \beta_\ell\ge0,\quad 0\le c_\ell<1.
\]
Split this sum as
\[
 f=\sum_{c_\ell<1/2}\beta_\ell\psi_{c_\ell},
 \qquad
 g=\sum_{c_\ell\ge1/2}\beta_\ell\psi_{c_\ell}.
\]
Then \(\widetilde h=f+g\), both functions are nondecreasing,
and \(g=0\) on \([0,1/2]\). We will use the shifted weight
comparison for \(f\) and the unshifted comparison for \(g\).

Let \(P,E_1,E_s\) be as in
Proposition~\ref{prop:weight-estimates}.
For \(j>P\), its pointwise bounds give
\[
 u_j\ge
 \begin{cases}
 f(y_j)+g(x_j),&j\notin E_s,\\
 f(x_{j-1})+g(x_j),&j\in E_s\setminus E_1,\\
 f(x_{j-1})+g(x_{j-1}),&j\in E_1.
 \end{cases}
\]
For \(j\le P\), we have \(x_j\le y_j<1/2\), so \(g\)
vanishes at both points and the total error bound becomes
\[
 \sum_{j=1}^P\bigl(f(y_j)-u_j\bigr)_+\le s f(1/2).
\]

For a function \(v\), put
\[
 \begin{aligned}
 Q_0[v]&=v(1)+\sum_{j=1}^a v(x_j)
                  -N\int_0^1v(\sqrt{x})\,dx,\\
 Q_s[v]&=v(1)+\sum_{j=1}^a v(y_j)
                  -N\int_0^1v(\sqrt{x})\,dx.
 \end{aligned}
\]
The first expression agrees with the notation of
Lemma~\ref{lem:hinge-quadrature}, since \(N=a+2\).
Summing the preceding bounds and using
\(u_{a+1}=1\ge\widetilde h(1)\), we obtain
\begin{equation}\label{eq:scalar-functional-lower-bound}
 \begin{aligned}
 W-N\int_0^1\widetilde h(\sqrt{x})\,dx
 \ge{}&Q_s[f]-s f(1/2)\\
 &-\sum_{j\in E_s}\bigl(f(y_j)-f(x_{j-1})\bigr)\\
 &+Q_0[g]-\sum_{j\in E_1}\bigl(g(x_j)-g(x_{j-1})\bigr).
 \end{aligned}
\end{equation}
By linearity, it remains to show that the contribution of
each summand \(\psi_c\) is positive.

The increments \(x_j-x_{j-1}\) decrease, so any \(M\)
distinct increments sum to at most \(\sqrt{M/a}\).
Since \(\psi_c\) is nondecreasing and \(1\)-Lipschitz,
\begin{equation}\label{eq:high-exception-loss}
 \sum_{j\in E_1}\bigl(\psi_c(x_j)-\psi_c(x_{j-1})\bigr)
 \le\min\{1-c,\sqrt{k/a}\}.
\end{equation}
The bound \(1-c\) follows by telescoping over all indices.
Each difference \(y_j-x_{j-1}\) contains at most \(s+1\)
consecutive increments. Summing separately at each offset gives
\begin{equation}\label{eq:low-exception-loss}
 \sum_{j\in E_s}\bigl(\psi_c(y_j)-\psi_c(x_{j-1})\bigr)
 \le(s+1)\sqrt{|E_s|/a}
 \le(s+1)^{3/2}\sqrt{k/a}\le s/8.
\end{equation}

For \(0\le c<1/2\), shifting the sum gives
\begin{equation}\label{eq:shifted-hinge-lower-bound}
 \begin{aligned}
 Q_s[\psi_c]
 &=Q_0[\psi_c]+s(1-c)-\sum_{j=1}^s\psi_c(x_j)\\
 &\ge Q_0[\psi_c]+s\bigl(1-c-\sqrt{s/a}\bigr).
 \end{aligned}
\end{equation}
Using \(\psi_c(1/2)=1/2-c\), \(a\ge16s\), and
\eqref{eq:low-exception-loss}, we obtain
\[
 \begin{aligned}
 Q_s[\psi_c]-s\psi_c(1/2)
 &-\sum_{j\in E_s}\bigl(\psi_c(y_j)-\psi_c(x_{j-1})\bigr)\\
 &\ge Q_0[\psi_c]+\frac{s}{2}-s\sqrt{s/a}-\frac{s}{8}\\
 &\ge Q_0[\psi_c]+\frac{s}{8}>0,
 \end{aligned}
\]
where the last inequality follows from
Lemma~\ref{lem:hinge-quadrature}.

For \(1/2\le c<1\), the same lemma and
\eqref{eq:high-exception-loss} give
\[
 Q_0[\psi_c]
 -\sum_{j\in E_1}\bigl(\psi_c(x_j)-\psi_c(x_{j-1})\bigr)
 >
 \begin{cases}
 1/12-\sqrt{k/a},&1/2\le c\le4/5,\\
 0,&4/5<c<1.
 \end{cases}
\]
The first bound is positive because
\eqref{eq:scalar-numerical-conditions} and \(s\ge48\) imply
\(k/a<1/144\).

Multiplying these inequalities by \(\beta_\ell\) and summing
makes the right-hand side of
\eqref{eq:scalar-functional-lower-bound} strictly positive.
Indeed, at least one coefficient is positive because
\(\widetilde h\ne0\). Therefore
\[
 W>N\int_0^1\widetilde h(\sqrt{x})\,dx
 \ge N\int_0^1h(\sqrt{x})\,dx,
\]
as required.
\end{proof}

Combining this proposition with
Proposition~\ref{prop:scalar-chow-criterion}, we obtain the following
numerical criterion.

\begin{corollary}\label{prop:numerical-stability-criterion}
Put \(s=16(2d+1)\), and suppose that \(k\) bounds the number
of subspaces with \(B^2>0\) and \(\kappa(B,V)\le2s+1\)
in every weight filtration. If \(r\ge3\) satisfies
\begin{equation}\label{eq:numerical-stability-conditions}
 \begin{gathered}
 r\ge32(s+1),\qquad
 \frac{s(1+\log(2dr))}{r}\le1,\qquad a=dr^2\ge16s,\\
 \frac{k}{a}\le
 \frac{s^2}{64(s+1)^3},
 \end{gathered}
\end{equation}
then \((X,L^{\otimes r})\) is Chow stable.
\end{corollary}

\subsection{An Ample Divisor on the Resolution}
\label{subsec:effective-auxiliary-divisor}

The subspace count depends on the degree of an auxiliary ample
divisor on the minimal resolution. To obtain a bound depending
only on \(d\), we choose such a divisor as follows.

\begin{lemma}\label{lem:effective-auxiliary-divisor}
Let \((X,L)\) be a complex polarized K3 surface with at worst
ADE singularities and \(L^2=2d\). Its minimal resolution
\(\pi:\widetilde X\to X\) admits an integral ample divisor \(H\) with
\[
 H\cdot\pi^*L=292d.
\]
\end{lemma}

\begin{proof}
Put \(\ell=\pi^*L\). The exceptional curves generate a direct
sum of negative-definite ADE lattices of total rank at most
\(19\). In each component, take the sum of the positive roots,
viewed as an effective exceptional divisor, and let \(E\) be
the sum of these divisors. The Weyl vector identities give
\[
 \ell\cdot E=0,\qquad E\cdot E_i=-2,\qquad
 M:=-E^2=\sum_\nu\frac{n_\nu h_\nu(h_\nu+1)}{3},
\]
where \(E_i\) runs through the exceptional curves, and
\(n_\nu,h_\nu\) are the rank and Coxeter number of each
component. The identities are recalled in
\cite[Remark~6.6 and Proposition~6.12]{Heckman}.
Since \(\sum_\nu n_\nu\le19\) and \(h_\nu\le36\), we have
\(M\le444\cdot19=8436\).

Set \(H=146\ell-E\). Then
\[
 H\cdot\ell=292d,\qquad H^2=2\cdot146^2d-M>0,
 \qquad H\cdot E_i=2.
\]
For any other irreducible curve \(C\), put
\(q=\ell\cdot C\ge1\). Adjunction gives \(C^2\ge-2\),
and the Hodge index theorem on \(\ell^\perp\) yields
\[
 (E\cdot C)^2
 \le M\left(\frac{q^2}{2d}-C^2\right)
 \le8436\left(\frac{q^2}{2}+2\right)
 \le21090q^2<146^2q^2.
\]
Thus \(H\cdot C>0\), and the Nakai--Moishezon criterion
shows that \(H\) is ample. For smooth \(X\), the given
polarization \(L\) already has the required ampleness.
\end{proof}

\subsection{An Effective Bound for the Exponent}
\label{sec:effective-degree-bound}
\label{sec:effective-scalar-bound}

We now substitute the explicit subspace count into
Corollary~\ref{prop:numerical-stability-criterion}, keeping track
of the degree of the auxiliary ample divisor.

\begin{proposition}\label{prop:effective-exponent}
Suppose that \(\widetilde X\) admits an integral ample divisor
\(H\) with \(H\cdot\pi^*L\le2dt\), where \(t\ge1\).
Put \(s=16(2d+1)\). Then \((X,L^{\otimes r})\) is Chow stable
for every integer
\begin{equation}\label{eq:effective-exponent-polynomial}
 r\ge3072t s^3=3\cdot2^{22}t(2d+1)^3.
\end{equation}
\end{proposition}

\begin{proof}
We verify the conditions of
Corollary~\ref{prop:numerical-stability-criterion}.
Since \(s\ge48\) and \(r\ge3072t s^3\), we have
\(r\ge32(s+1)\) and \(a=dr^2\ge16s\).
The function \((1+\log(2dr))/r\) is decreasing for \(2dr>1\).
Using \(r\ge s^3\) and \(2d\le s\), we obtain
\[
 \frac{s(1+\log(2dr))}{r}
 \le\frac{1+4\log s}{s^2}\le1.
\]

It remains to bound the number of exceptional subspaces.
Apply \eqref{eq:explicit-subspace-count} with \(\eta=2s+1\).
Since \(r\ge s^3\), the threshold is \(R=16dr\), and
\[
 \frac{9\eta^2+\eta}{2}+1
 =18s^2+19s+6\le19s^2.
\]
Together with \(c=H\cdot\pi^*L\le2dt\), this gives
\[
 k\le8dr+2+19s^2(2dtr+1)\le40dtr s^2.
\]
Consequently,
\[
 \frac{64(s+1)^3}{s^2}\frac{k}{a}
 \le\frac{2560t(s+1)^3}{r}
 \le\frac56\left(1+\frac1s\right)^3<1.
\]
This verifies the last condition and proves the assertion.
\end{proof}

\subsection{Proof of Theorem~\ref{thm:main}}
\label{sec:proof-main}
For an ADE K3 surface,
Lemma~\ref{lem:effective-auxiliary-divisor} supplies \(H\) with
\(H\cdot\pi^*L=2d\cdot146\).
Proposition~\ref{prop:effective-exponent}, applied with
\(t=146\), therefore gives Chow stability whenever
\[
 r\ge2^{31}(2d+1)^3,
\]
since \(3\cdot146\cdot2^{22}<2^{31}\).
If \(X\) is smooth, applying the same proposition with
\(H=L\) and \(t=1\) gives Chow stability whenever
\[
 r\ge2^{24}(2d+1)^3.
\]
This proves Theorem~\ref{thm:main}.
\qed

\bibliographystyle{alpha}
\bibliography{main}

\end{document}